\documentclass{amsart}
\usepackage{graphicx}

\newtheorem{thm}{Theorem}[section]
\newtheorem{cor}[thm]{Corollary}

\theoremstyle{definition}

\theoremstyle{remark}

\numberwithin{equation}{section}
\begin{document}

\title[Multiserver queueing systems with retrials and losses]
{Multiserver queueing systems with retrials and losses}%
\author{Vyacheslav M. Abramov}%
\address{School of Mathematical Sciences, Monash University, Building 28M,
Clayton, Victoria 3800, Australia}%
\email{vyacheslav.abramov@sci.monash.edu.au}%

%\thanks{}%

%\date{}%
%\dedicatory{}%
%\commby{}%
% ----------------------------------------------------------------
\begin{abstract}
The interest to retrial queueing systems is due to their
application to telephone systems. The paper studies multiserver
retrial queueing systems with $n$ servers. Arrival process is a
quite general point process. An arriving customer occupies one of
free servers. If upon arrival all servers are busy, then the
customer waits for his service in orbit, and after random time
retries more and more to occupy a server. The orbit has one
waiting space only, and arriving customer, who finds all servers
busy and the waiting space occupied, losses from the system. Time
intervals between possible retrials are assumed to have arbitrary
distribution (the retrial scheme is exactly explained in the
paper). The paper provides analysis of this system. Specifically
the paper studies optimal number of servers to decrease the loss
proportion to a given value. The representation obtained for loss
proportion enables us to solve the problem numerically. The
algorithm for numerical solution includes effective simulation,
which meets the challenge of rare events problem in simulation.
\end{abstract}
\maketitle
\tableofcontents
% ----------------------------------------------------------------
\section{Introduction}\label{Introduction}
The paper studies multiserver queueing systems with retrials and
losses. Retrial queueing systems are an object of investigation in
numerous papers (e.g. see reviews Artalejo \cite{Artalejo 1999},
Artalejo and Falin \cite{Artalejo and Falin 2002}, Falin
\cite{Falin 1990}). The interest to these systems is due to their
application to telephone systems.

During the last years there has been an especial interest to
queueing systems with retrial and losses. The different queueing
systems of this type has been studied in \cite{Bocharov et al
2001}, \cite{Kovalenko 2002}, \cite{Mandelbaum et al 2002}, and
other papers.

The most of known studies in this direction are devoted to
analytical solution of one or other problem. Bocharov et al
\cite{Bocharov et al 2001} studied multiclass single-server M/G/1
queueing system with finite buffer and limited number of places in
orbit. They deduced relationship between steady state
distributions of this system and the similar system with only one
class of customers. Kovalenko \cite{Kovalenko 2002} studied the
loss probability in M/G/$n$ retrials queueing system with two
customer classes. Mandelbaum et al \cite{Mandelbaum et al 2002}
studied queueing system with retrials and losses with
time-dependent parameters in framework of Markovian service
network (see \cite{Mandelbaum Massey and Reiman 1998}. They
provide fluid and diffusion approximations to obtain numerical
characteristics for queue-length and virtual waiting time
processes. Atencia and Phong \cite{Atencia and Phong 2004} and
Bocharov, Phong and Atencia \cite{Bocharov Phong and Atencia 2001}
also studied characteristics of queueing systems with retrial and
losses with the aid of analytic methods.

However explicit representations can be deduced for only
restricted classes of retrial queueing systems, and application of
these results for real systems is very difficult.

Recently Abramov \cite{Abramov 2006} studied the main multiserver
retrial queueing system, where input stream was a quite general
point process (see \cite{Abramov 2006} for more details). The
martingale approach, that has been used there, can be adapted to
systems with retrials and losses as well, and the present paper is
just devoted to analysis of these systems. The exact description
of the system is as follows.

\smallskip
$\bullet$ The arrival process $A(t)$ is a point process (the
assumptions on this process are given later).

\smallskip
$\bullet$ There are $n$ servers, and an arriving customer occupies
one of free servers.

\smallskip
$\bullet$ If upon arrival all servers are busy, but the secondary
queue, orbit, having only one space is free, then the customer
 occupies the orbit and retries more and more to occupy a server.

\smallskip
$\bullet$ A customer, who upon arrival finds all servers busy and
the orbit occupied, loses from the system without delay.

\smallskip
$\bullet$ A service time of each customer is exponentially
distributed random variable with parameter $\mu$.

\smallskip
$\bullet$ A time between retrials is arbitrary distributed. It is
generated by point process $D(t)$ (the details on this assumption
are given later). The processes $D(t)$ and $A(t)$ are assumed to
be disjoint (that is the probability of simultaneous jump of these
processes is equal to 0).

\smallskip
Notice, the significant difference between the main multiserver
model of \cite{Abramov 2006} and the model considered in the
present paper is also that retrial times in the present model are
generally distributed, while in \cite{Abramov 2006} retrial times
of each customer in orbit are exponentially distributed.

There is also difference in the assumptions. In the case of the
main multiserver retrial queueing systems considered in
\cite{Abramov 2006} the point process $A(t)$ must satisfy the
condition
\begin{equation}
\label{I1}\mathbb{P}\left\{\lim_{t\to\infty}\frac{A(t)}{t}=\lambda<n\mu\right\}
=1,
\end{equation}
leading to the stability of the system.

In the case of the system considered in the present paper the
number of customers in the systems is always bounded, and the
stability condition \eqref{I1} is not longer required, and the
following more weaker conditions of convergence are used
\begin{equation}
\label{I2}\mathbb{P}\left\{\lim_{t\to\infty}\frac{A(t)}{t}=\lambda\right\}=1,
\end{equation}
and
\begin{equation}
\label{I3}\mathbb{P}\left\{\lim_{t\to\infty}\frac{D(t)}{t}=\Delta\right\}=1,
\end{equation}
i.e. $\lambda$, $\Delta$, $\mu$ and $n$ are not  allied to one
another.

Conditions \eqref{I2} and \eqref{I3} are technical conditions,
that are then used by the methods of our analysis. (We often use
the Lebesgue theorem on dominated convergence.)

The main aim of this paper is the numerical solution of a concrete
problem associated with specific performance measure. The main
question that we want to answer in this paper is \textit{under
what number of servers the loss proportion will be less than a
given small value?} The system that will be studied in this paper
is described by quite general point processes, and explicit
analytic solution for this system in general case is impossible.
At the same time the direct simulation of this system in order to
answer the aforementioned question is not available either, since
as the loss probability is chosen very small, we should deal with
simulating rare event problem. The simulating rare events problem
is a well-known problem having notable attention in the literature
(e.g. \cite{Heidelberger 1995}, \cite{Shahabuddin 1995} and many
others). Then the main result of this paper is an answer the
question \textit{how to simulate this system?}

The main results of the present paper are represented by relation
\eqref{SDN8} and Theorem \ref{thm1}. These results enable us to
\textit{effectively simulate} models to study their most
significant performance characteristic, the loss proportion. By
effective simulation we mean such an approach that make the
challenge of the rare events problem in simulation.

 In our case, we discuss three relations for loss proportions:
 \eqref{SDN8}, \eqref{SDN9} and \eqref{SDN10}. Relations \eqref{SDN9} and
 \eqref{SDN10} are not effective.  For example,
\eqref{SDN10}
 requires to simulate the events $\{Q_1(s-)=n\}$ and
 $\{Q_2(s-)=1\}$ directly. These events become rare as $n$ increases indefinitely.
 In contrast relation \eqref{SDN8} is effective. Its expression
 contains $\lim_{t\to\infty}t^{-1}\mathbb{E}\int_0^t Q(s)\mbox{d}s$, which
 according to Theorem \ref{thm1} can be easily approximated, so
 that the probabilities of rare events presented in calculation of this expression
 can be removed. Specifically, there can be used a change of the Poisson
 measure in \eqref{SDN2}
 for calculation of \eqref{SDN8}, i.e. replace one Poisson measure to another one and
 use the Radon-Nikodim derivative. In other words, the Poisson process
 under which the required characteristic of the system has a small
 probability, is replaced by other Poisson process under which
 the usual simulation becomes available, and then the required small
 probability is recalculated by using the Radon-Nikodim derivative
 (see e.g. \cite{Rubinstein and Shapiro 1993} or \cite{Melamed and Rubinstein
 1998}).
 Recall that increments of a Poisson process are exponentially distributed,
and for two exponentially distributed random variables with
parameters $\vartheta_1$ and $\vartheta_2$, the Radon-Nikodim
derivative is
$\frac{\mbox{d}\mathbb{P}_{\vartheta_1}(x)}{\mbox{d}\mathbb{P}_{\vartheta_2}(x)}
=\frac{\vartheta_1}{\vartheta_2}
\mbox{e}^{(\vartheta_2-\vartheta_1)x}$, where
$\mathbb{P}_{\vartheta_k}(x)$, $k=1,2$, are corresponding
exponential probability distributions.

 Relations \eqref{SDN9} and \eqref{SDN10} are derivative from
\eqref{BE2} and \eqref{BE3} correspondingly.

The presence of these equations is illustrative in order to show
the difference between direct simulation and simulation based on
\eqref{SDN8} and Theorem \ref{thm1}.

 The paper is organized as follows. In Section \ref{Basic
 equations} we discuss basic equations giving us then three
 aforementioned relations \eqref{SDN8}, \eqref{SDN9} and \eqref{SDN10} for loss proportions.
In Section \ref{Normalization} we use semimartingale decomposition
and Lenglart-Rebolledo inequality to derive then crucial relation
\eqref{SDN8}. In Section \ref{Frequences} we prove Theorem
\ref{thm1} and the corollary from this theorem related to
particular case where the point processes $A(t)$ and $D(t)$ are
Poisson. Specifically, Theorem \ref{thm1} establishes
relationships for the proportions
$\lim_{t\to\infty}t^{-1}\int_0^t\mathbb{P}\{Q_1(s)=i,
Q_2(s)=j\}\mbox{d}s$ ($i=0,1,\ldots,n-1$; $j=0,1$). Section
\ref{Performance} is devoted to analysis of the performance
measure based on loss proportion. There is formulated the problem
of optimization of the proposed performance measure and described
the algorithms for its solution. Section \ref{Numerical examples}
provides two numerical examples.  The concluding remarks are in
Section \ref{Concluding remarks}.

\section{Basic equations}\label{Basic equations}
All point processes considered in this paper are assumed to be
right-continuous having the left-side limits.

The following notation is used. The arrival process is denoted
$A(t)$. The queueing process (the number of occupied servers
process) is denoted $Q_1(t)$. The orbit state is denoted $Q_2(t)$.
It can take only the values 0 and 1 correspondingly to the cases
of empty and occupied orbit space. The loss process is denoted
$Q_3(t)$. Specifically, $Q_3(t)$ is the cumulated number of losses
up to time $t$. The moments of retrials are governed by point
process $D(t)$, that is, the retrial moments are $t_1$, $t_2$,
\ldots and $D(t)$ = $\min\{k: t_k\leq t\}$. The sequence of
independent Poisson processes all with the same rate $\mu$ is
denoted $\pi_l(t)$.

Then we have the following equations ($\mathbf{I}\{\mathcal{A}\}$
denotes the indicator of event $\mathcal{A}$):
\begin{eqnarray}
\label{BE1}
Q_1(t)+Q_2(t)+Q_3(t)&=&A(t)-\int_0^t\sum_{l=1}^{n}\mathbf{I}\{Q_1(s-)\geq
l\}\mbox{d}\pi_l(s),\\
\label{BE2}Q_2(t)+Q_3(t)&=&\int_0^t\mathbf{I}\{Q_1(s-)={n}\}\mbox{d}A(s)\\
&&-
\int_0^t\mathbf{I}\{Q_1(s-)\neq
n\}\mathbf{I}\{Q_2(s-)=1\}\mbox{d}D(s),\nonumber\\
\label{BE3}Q_3(t)&=&\int_0^t\mathbf{I}\{Q_1(s-)={n}\}\mathbf{I}\{Q_2(s-)=1\}\mbox{d}A(s).
\end{eqnarray}
Clearly, that the meaning of the left-hand side of relation
\eqref{BE1} is the total number of customers in the system
including customers in service, one in orbit (in the case of
occupied place) as well as all losses up to time $t$. Notice, that
the right-hand side of \eqref{BE1} is as the right-hand side of
(2.1) in \cite{Abramov 2006}. The first term of the right-hand
side of \eqref{BE2} means the total number of arrivals to busy
system where all servers are occupied. The second term of the
right-hand side of \eqref{BE2} is the subtracted number of
successfully retrial customers up to time $t$ from the total
number of arrivals to busy system during the same time $t$.
Successful retrials occur in the cases where immediately before
retrial instants there is at least one busy server. Notice, that
relation \eqref{BE2} is structured similarly to corresponding
relation (2.2) of \cite{Abramov 2006}. Relation \eqref{BE3}
characterizes the number of losses. Losses occur in the cases
where immediately before arrival all servers are busy and the
place in orbit is occupied. This is indicated by the right-hand
side of \eqref{BE3}.

\section{Semimartingale decompositions and normalization}\label{Normalization}
In this section we discuss the normalized processes by dividing
these processes by $t$. As $t$ increases to infinity, the
behaviour of the process $Q_3(t)$ divided by $t$ is the
significant performance characteristics characterizing losses per
time unit.

For the purpose of study normalized characteristics we need in
semimartingale decomposition of the processes (e.g. Liptser and
Shiryayev \cite{Liptser and Shiryayev 1977-1978}, Jacod and
Shiryayev \cite{Jacod and Shiryayev 1987}).

The processes $A(t)$, $D(t)$ and $\pi_l(t)$, $l=1,2,\ldots,m$ all
are assumed to be given on the common filtered probability space
($\Omega$, $\mathcal{F}$, $\mathbf{F}=(\mathcal{F}_t)_{t\geq 0}$,
$\mathbb{P}$). All these processes are semimartingales. The
semimartingale decomposition for arbitrary point process $X(t)$ is
written $X(t)=\widehat X(t)+M_X(t)$, where $\widehat X(t)$ is a
compensator and $M_X(t)$ is a local square integrable martingale.
For example: $A(t)$=$\widehat A(t)$+$M_A(t)$. In some cases for
local square integrable  martingales we use the notation $M$ with
indexes such as $M_{i,j,k}(t)$. Such types of notation will be
especially explained.

Observing first \eqref{BE1}, rewrite it
\begin{eqnarray}
\label{SDN1} Q_1(t)+Q_2(t)+Q_3(t)&=&A(t)-C(t),\\
\label{SDN2}C(t)&=&\int_0^t\sum_{l=1}^{n} \textbf{I}\{Q_1(s-)\geq
l\}\mbox{d}\pi_l(s).
\end{eqnarray}
The semimartingale $C(t)$ admits Doob-Meyer decomposition
\begin{equation}
\label{SDN3} C(t)=\widehat C(t)+M_C(t).
\end{equation}
In turn, the compensator $\widehat C(t)$ admits the representation
\begin{equation}
\label{SDN4} \widehat C(t)=\mu\int_0^t Q_1(s)\mbox{d}s
\end{equation}
(for details see Dellacherie \cite{Dellacherie 1972}, Liptser and
Shiryayev \cite{Liptser and Shiryayev 1977-1978}, \cite{Liptser
and Shiryayev 1989}, Theorem 1.6.1). Therefore, in view of
\eqref{SDN2}, \eqref{SDN3} and \eqref{SDN4}, one can rewrite
\eqref{SDN1} as follows:
\begin{equation}
\label{SDN5}Q_1(t)+Q_2(t)+Q_3(t)=A(t)-\mu\int_0^t
Q_1(s)\mbox{d}s-M_C(t).
\end{equation}
Let us pass to normalized processes. For arbitrary process $X(t)$,
for $t>0$ its normalization is denoted by small letter, so
\begin{equation*}
x(t)=\frac{X(t)}{t}.
\end{equation*}
Therefore, passing to normalized processes in \eqref{SDN5} can we
written as
\begin{equation}
\label{SDN6}q_1(t)+q_2(t)+q_3(t)=a(t)-\frac{\mu}{t}\int_0^t
Q_1(s)\mbox{d}s-m_C(t).
\end{equation}
Assume now that $t$ increases to infinity. Then, because of
stochastic boundedness of $Q_1(t)$ and $Q_2(t)$ the terms $q_1(t)$
and $q_2(t)$ vanish with probability 1, and the left-hand side of
\eqref{SDN6} in limit looks
$\mathbb{P}^{\_}\lim_{t\to\infty}q_3(t)$.

Next, according to \eqref{I2}, as $t$ increases to infinity,
$a(t)$ converges with probability 1 to $\lambda$. The term
$m_C(t)$ vanishes in probability, since according to
Lenglart-Rebolledo inequality we have:
\begin{equation}\label{SDN7}
\begin{aligned}
\mathbb{P}\{|m_C(t)|>\delta\}&\leq\mathbb{P}\left\{\sup_{0<s\leq
t}\left|\frac{sm_C(s)}{t}\right|>\delta\right\}\\
&=\mathbb{P}\left\{\sup_{0<s\leq t}\left|M_C(s)\right|>\delta
t\right\}\\
&\leq\frac{\epsilon}{\delta^2}+\mathbb{P}\left\{\sum_{l=1}^{n}\pi_l(t)>\epsilon
t^2\right\}\\
&=\frac{\epsilon}{\delta^2}+\mathbb{P}\left\{\frac{1}{t}\sum_{l=1}^{n}\pi_l(t)>\epsilon
t\right\},
\end{aligned}
\end{equation}
and because of arbitrariness of $\epsilon$ in \eqref{SDN7} the
statement is right.

Therefore, applying the Lebesgue theorem on dominated convergence
\begin{equation}
\label{SDN8}
\begin{aligned}
\mathbb{P}^{\_}\lim_{t\to\infty}q_3(t)&=\lambda-\mu\left(\mathbb{P}^{\_}\lim_{t\to\infty}
\frac{1}{t}\int_0^tQ_1(s)\mbox{d}s\right)\\
&=\lambda-\mu\left(\lim_{t\to\infty}
\frac{1}{t}\mathbb{E}\int_0^tQ_1(s)\mbox{d}s\right)
\end{aligned}
\end{equation}
The meaning of the term of the right-hand side of \eqref{SDN8},
taken in brackets, is the expected number of occupied servers
during the long-run history, or expected stationary number of
occupied servers.

Relations \eqref{BE2} and \eqref{BE3} provide additional
information on the behaviour of the process $q_3(t)$.
Specifically, we have
\begin{equation}
\label{SDN9}
\begin{aligned}
\mathbb{P}^{\_}\lim_{t\to\infty}q_3(t)=&\lim_{t\to\infty}
\frac{1}{t}\mathbb{E}\int_0^t\mathbf{I}\{Q_1(s-)={n}\}\mbox{d}A(s)\\
&-\lim_{t\to\infty}
\frac{1}{t}\mathbb{E}\int_0^t\mathbf{I}\{Q_1(s-)\neq
{n}\}\mathbf{I}\{Q_2(s-)=1\}\mbox{d}D(s)
\end{aligned}
\end{equation}
and
\begin{equation}
\label{SDN10}\mathbb{P}^{\_}\lim_{t\to\infty}q_3(t)=
\lim_{t\to\infty} \frac{1}{t}
\mathbb{E}\int_0^t\mathbf{I}\{Q_1(s-)={n}\}\mathbf{I}\{Q_2(s-)=1\}\mbox{d}A(s).
\end{equation}

In the particular case where the process $A(t)$ is Poisson from
\eqref{SDN10} we obtain:
\begin{equation}
\label{SDN11}
\begin{aligned}
\mathbb{P}^{\_}\lim_{t\to\infty}q_3(t)&= \lambda \lim_{t\to\infty}
\frac{1}{t}\mathbb{E}\int_0^t\mathbf{I}\{Q_1(s)={n}\}\mathbf{I}\{Q_2(s)=1\}\mbox{d}s\\
&=\lambda \lim_{t\to\infty}
\frac{1}{t}\int_0^t\mathbb{P}\{Q_1(s)={n}, Q_2(s)=1\}\mbox{d}s,
\end{aligned}
\end{equation}
and from \eqref{SDN9} we obtain
\begin{equation}
\label{SDN12}
\begin{aligned}
\mathbb{P}^{\_}\lim_{t\to\infty}q_3(t)=&\lambda\lim_{t\to\infty}
\frac{1}{t}\int_0^t\mathbb{P}\{Q_1(s)={n}\}\mbox{d}s\\
&-\lim_{t\to\infty}
\frac{1}{t}\mathbb{E}\int_0^t\mathbf{I}\{Q_1(s-)\neq
{n}\}\mathbf{I}\{Q_2(s-)=1\}\mbox{d}D(s).
\end{aligned}
\end{equation}
In the case where both $A(t)$ and $D(t)$ are Poisson, $q_3(t)$
converges to its limit with probability 1 as $t$ increases to
infinity, because both $m_A(t)$ and $m_D(t)$ vanish with
probability 1. Specifically, in this case we obtain:
\begin{equation}
\begin{aligned}
\label{SDN13}&\mathbb{P}\left\{\lim_{t\to\infty}q_3(t)=
\lambda\lim_{t\to\infty}
\frac{1}{t}\int_0^t\mathbb{P}\{Q_1(s)={n}\}\mbox{d}s\right.\\
&\ \ \ \ -\Delta\lim_{t\to\infty}
\left.\frac{1}{t}\int_0^t\mathbb{P}
\{Q_1(s)\neq {n}, Q_2(s)=1\}\mbox{d}s\right\}\\
&=\mathbb{P}\left\{\lim_{t\to\infty}q_3(t)=
\lambda\lim_{t\to\infty}
\mathbb{P}\{Q_1(t)={n}\}\right.\\
&\ \ \ \ -\Delta\left.\lim_{t\to\infty} \mathbb{P}
\{Q_1(t)\neq {n}, Q_2(t)=1\}\right\}\\
&=1.
\end{aligned}
\end{equation}
There is use of the fact of existence of the limiting stationary
probabilities in \eqref{SDN13} as $t\to\infty$, as well as the
equalities
\begin{eqnarray*}
\lim_{t\to\infty} \mathbb{P}\{Q_1(t)={n}\}&=&\lim_{t\to\infty}
\frac{1}{t}\int_0^t\mathbb{P}\{Q_1(s)={n}\}\mbox{d}s,\\
\lim_{t\to\infty} \mathbb{P} \{Q_1(t)\neq {n}, Q_2(t)=1\}&=&
\lim_{t\to\infty} \frac{1}{t}\int_0^t\mathbb{P} \{Q_1(s)\neq {n},
Q_2(s)=1\}\mbox{d}s.
\end{eqnarray*}

\section{Limiting frequencies and distributions}\label{Frequences}
The results obtained in the previous section do not provide
completed information about losses in this system. We have no
information about the quantity $\lim_{t\to\infty}t^{-1}\int_0^t
Q_1(s)\mbox{d}s$ entering to the right-hand side of \eqref{SDN8}.
The information about the right-hand sides of \eqref{SDN9} and
\eqref{SDN10} is unknown as well.

Therefore in this section we derive equations for the following
frequencies:
\begin{equation}\label{LDF1}
\begin{aligned}
&\lim_{t\to\infty}\frac{1}{t}\int_0^t\mathbb{P}\{Q_1(s)=i,
Q_2(s)=j\}\mbox{d}s,\\
&i=0,1,\ldots,{n}; \ j=0,1.
\end{aligned}
\end{equation}
For this purpose we introduce the processes:
\begin{equation}\label{LDF2}
\begin{aligned}
&I_{i,j}(t)=\mathbf{I}\{Q_1(t)=i \ \cap \ Q_2(t)=j\},\\
&i=0,1,\ldots,{n}; \ j=0,1,
\end{aligned}
\end{equation}
taking the value 0 if at least one of the indexes $i$ or $j$ is
negative.

The jump of a point process is denoted by adding $\triangle$
\footnote{This is the triangle sign. We hope that this sign will
not be mixed with the parameter (capital Greek letter)  $\Delta$}.
For example, $\triangle A(t)$ is a jump of $A(t)$, $\triangle
\pi_l(t)$ is a jump for $\pi_l(t)$ and so on.

Denote
\begin{equation*}
\Pi_i(t)=\sum_{l=1}^i\pi_l(t).
\end{equation*}

Then we have the following equations:
\begin{eqnarray}
\label{LDF3}&&\mathbf{I}\{Q_1(t-)+\triangle Q_1(t)=i \ \cap
\ Q_2(t)=Q_2(t-)=j\}\\
&=&I_{i-1,j}(t-)\triangle
A(t)+I_{i+1,j}(t-)\triangle\Pi_{i+1}(t)\nonumber\\
&&+\ I_{i,j}(t-)[1-\triangle
A(t)][1-\triangle\Pi_i(t)][1-j\triangle D(t)],\nonumber\\
&& \ \ \ \ i=0,1,\ldots,{n}-1; \ j=0,1,\nonumber
\end{eqnarray}
\begin{equation}
\label{LDF4}
\begin{aligned}
&\mathbf{I}\{Q_1(t-)+\triangle Q_1(t)=i \ \cap Q_2(t-)=1 \ \cap
\ Q_2(t)=0\}\\
&=I_{i-1,1}(t-)\triangle D(t),\\
& \ \ \ \ i=0,1,\ldots,{n}-1,
\end{aligned}
\end{equation}
\begin{equation}
\label{LDF5} \begin{aligned} &\mathbf{I}\{Q_1(t-)+\triangle
Q_1(t)={n} \ \cap
\ Q_2(t)=Q_2(t-)=0 \}\\
&=I_{{n}-1,0}(t-)\triangle
A(t)+I_{{n},0}(t-)[1-\triangle\Pi_{{n}}(t)][1-\triangle A(t)],\\
\end{aligned}
\end{equation}
\begin{equation}
\label{LDF6}
\begin{aligned}
&\mathbf{I}\{Q_1(t-)+\triangle Q_1(t)={n} \ \cap Q_2(t-)=0 \ \cap
\
Q_2(t)=1\}\\
&=I_{{n},0}(t-)\triangle A(t),\\
\end{aligned}
\end{equation}
\begin{equation}
\label{LDF7}
\begin{aligned}
&\mathbf{I}\{Q_1(t-)+\triangle Q_1(t)={n} \ \cap
\ Q_2(t)=Q_2(t-)=1\}\\
&=I_{{n}-1,1}(t-)\triangle A(t)+I_{{n},1}(t-)\triangle
A(t)+I_{{n},1}(t-)[1-\triangle A(t)][1-\triangle \Pi_{n}],
\end{aligned}
\end{equation}
\begin{equation}
\label{LDF8}
\begin{aligned}
&\mathbf{I}\{Q_1(t-)+\triangle Q_1(t)={n} \ \cap
\ Q_2(t-)=1 \ \cap \ Q_2(t)=0\}\\
&=I_{{n}-1,1}(t-)\triangle D(t).
\end{aligned}
\end{equation}

Then,
\begin{equation}
\label{LDF9}
\begin{aligned}
\triangle I_{i,j}(t)=&\mathbf{I}\{Q_1(t-)+\triangle Q_1(t)=i \
\cap
\ Q_2(t)=Q_2(t-)=j\}\\
&+\mathbf{I}\{Q_1(t-)+\triangle Q_1(t)=i \ \cap \ Q_2(t-)\neq
Q_2(t)=j\}\\
&-I_{i,j}(t-),\\
& \ \ \ \ i=0,1,\ldots,{n}; \ j=0,1.
\end{aligned}
\end{equation}

Since
\begin{equation*}
\sum_{s\leq t}\triangle I_{i,j}=I_{i,j}(t)-I_{i,j}(0),
\end{equation*}
then taking into account \eqref{LDF3} - \eqref{LDF9} we have the
following.

For $i=0,1,\ldots,{n}-1$,
\begin{equation}
\label{LDF10}
\begin{aligned}
I_{i,j}(t)=&I_{i,j}(0)+\int_0^t[I_{i-1,j}(s-)-I_{i,j}(s-)]\mbox{d}A(s)\\
&-\int_0^tI_{i,j}(s-)\mbox{d}\Pi_i(s)-j\int_0^tI_{i,j}(s-)\mbox{d}D(s)\\
&+\int_0^tI_{i+1,j}(s-)\mbox{d}\Pi_{i+1}(s)+(1-j)\int_0^tI_{i-1,1}(s-)\mbox{d}D(s).
\end{aligned}
\end{equation}

For $i={n}$ we have the following pair of equations. In the case
$j=0$ we have:
\begin{equation}
\label{LDF11}
\begin{aligned}
I_{{n},0}(t)=&I_{{n},0}(0)+\int_0^t[I_{{n}-1,0}(s-)-I_{{n},0}(s-)]\mbox{d}A(s)\\
&-\int_0^tI_{{n},0}(s-)\mbox{d}\Pi_{n}(s)\\
&+\int_0^tI_{{n}-1,1}(s-)\mbox{d}D(s).
\end{aligned}
\end{equation}

In turn, in the case $j=1$ we have:
\begin{equation}
\label{LDF12}
\begin{aligned}
I_{{n},1}(t)=&I_{{n},1}(0)+\int_0^t I_{{n}-1,1}(s-)\mbox{d}A(s)\\
&+\int_0^t
I_{{n},0}(s-)\mbox{d}A(s)-\int_0^tI_{{n},1}(s-)\mbox{d}\Pi_{n}(s).
\end{aligned}
\end{equation}

By semimartingale decomposition, from \eqref{LDF10}, \eqref{LDF11}
and \eqref{LDF12} we obtain the following equations. In the cases
$i=0,1,\ldots,{n}-1$; $j=0,1$ we have:
\begin{equation}
\label{LDF13}
\begin{aligned}
I_{i,j}(t)=&I_{i,j}(0)+\int_0^t[I_{i-1,j}(s-)-I_{i,j}(s-)]\mbox{d}A(s)\\
&-i\mu\int_0^t I_{i,j}(s)\mbox{d}s+(i+1)\mu\int_0^t
I_{i+1,j}(s)\mbox{d}s\\
&-j\int_0^tI_{i,1}(s-)\mbox{d}D(s)+(1-j)\int_0^t
I_{i-1,1}(s-)\mbox{d}D(s)\\
&+M_{i,j}(t),
\end{aligned}
\end{equation}
with the martingale
\begin{equation}\label{LDF13A}
M_{i,j}(t)=-\int_0^tI_{i,j}(s-)\mbox{d}M_{\Pi_i}(s)+\int_0^tI_{i+1,j}(s-)\mbox{d}
M_{\Pi_{i+1}}(s).
\end{equation}
For $i={n}$ we have the following pair of equations. In the case
$j=0$ we have:
\begin{equation}
\label{LDF14}
\begin{aligned}
I_{{n},0}(t)=&I_{{n},0}(0)+\int_0^t[I_{{n}-1,0}(s-)-I_{{n},0}(s-)]\mbox{d}A(s)\\
&-{n}\mu\int_0^t I_{{n},0}(s)\mbox{d}s+\int_0^t
I_{{n}-1,1}(s-)\mbox{d}D(s)\\
&+M_{{n},0}(t),
\end{aligned}
\end{equation}
with the martingale
\begin{equation}\label{LDF14A}
M_{{n},0}(t)=-\int_0^tI_{{n},0}(s-)\mbox{d}M_{\Pi_{n}}(s).
\end{equation}
In turn, in the case $j=1$ we have:
\begin{equation}
\label{LDF15}
\begin{aligned}
I_{{n},1}(t)=&I_{{n},1}(0)+\int_0^t
I_{{n}-1,1}(s-)\mbox{d}A(s)-{n}\mu\int_0^t
I_{{n},1}(s)\mbox{d}s\\&+ \int_0^t
I_{{n},0}(s-)\mbox{d}A(s)+M_{{n},1}(t),
\end{aligned}
\end{equation}
with the martingale
\begin{equation}\label{LDF15A}
M_{{n},1}(t)=-\int_0^tI_{{n},1}(s-)\mbox{d}M_{\Pi_{n}}(s).
\end{equation}

\begin{thm}\label{thm1}
For the limiting state frequencies of the processes $Q_1(t)$ and
$Q_2(t)$ we have the following equations.

In the boundary cases $i=0$ for $j=0,1$ we have:
\begin{equation}
\label{LDF16}
\begin{aligned}
 &\mu\lim_{t\to\infty}\frac{1}{t}\int_0^t\mathbb{P}\{
Q_1(s)=1, Q_2(s)=j\}\mbox{d}s\\
&=\lim_{t\to\infty}\frac{1}{t}\mathbb{E}\int_0^t\mathbf{I}\{Q_1(s-)=0,
Q_2(s-)=j\}\mbox{d}A(s)\\
&\ \ \ \ +j \
\lim_{t\to\infty}\frac{1}{t}\mathbb{E}\int_0^t\mathbf{I}\{Q_1(s-)=0,
Q_2(s-)=j\}\mbox{d}D(s).
\end{aligned}
\end{equation}

In the cases $i=1,2,\ldots,{n}-1$, $j=0,1$ we have:
\begin{equation}
\label{LDF17}
\begin{aligned}
&i\mu\lim_{t\to\infty}\frac{1}{t}\int_0^t\mathbb{P}\{Q_1(s)=i,
Q_2(s)=j\}\mbox{d}s\\&\ \ \ \
-(i+1)\mu\lim_{t\to\infty}\frac{1}{t}\int_0^t\mathbb{P}\{Q_1(s)=i+1,
Q_2(s)=j\}
\mbox{d}s\\
&=\lim_{t\to\infty}\frac{1}{t}\mathbb{E}\int_0^t[\mathbf{I}\{Q_1(s-)=i-1,
Q_2(s-)=j\}\\
&\ \ \ \ \ \ \ \ \ \ \ \ \ \ \ -\mathbf{I}\{Q_1(s-)=i, Q_2(s-)=j\}]\mbox{d}A(s)\\
&\ \ \ \ +
(1-j)\lim_{t\to\infty}\frac{1}{t}\mathbb{E}\int_0^t\mathbf{I}\{Q_1(s-)=i-1,
Q_2(s-)=1\}\mbox{d}D(s)\\
&\ \ \ \ -
j\lim_{t\to\infty}\frac{1}{t}\mathbb{E}\int_0^t\mathbf{I}\{Q_1(s-)=i,
Q_2(s-)=1\}\mbox{d}D(s).
\end{aligned}
\end{equation}

In the case $i={n}$ and $j=0$ we have:
\begin{equation}
\label{LDF18}
\begin{aligned}
&{n}\mu\lim_{t\to\infty}\frac{1}{t}\int_0^t\mathbb{P}\{Q_1(s)={n},
Q_2(s)=0\}\mbox{d}s\\
&=\lim_{t\to\infty}\frac{1}{t}\mathbb{E}\int_0^t[\mathbf{I}\{Q_1(s-)={n}-1,
Q_2(s-)=0\}\\
&\ \ \ \ \ \ \ \ \ \ \ \ \ \ \ -\mathbf{I}\{Q_1(s-)={n}, Q_2(s-)=0\}]\mbox{d}A(s)\\
&\ \ \ \ +
\lim_{t\to\infty}\frac{1}{t}\mathbb{E}\int_0^t\mathbf{I}\{Q_1(s-)={n},
Q_2(s-)=1\}\mbox{d}D(s).
\end{aligned}
\end{equation}

In the case $i={n}$ and $j=1$ we have:
\begin{equation}
\label{LDF19}
\begin{aligned}
&{n}\mu\lim_{t\to\infty}\frac{1}{t}\int_0^t\mathbb{P}\{Q_1(s)={n},
Q_2(s)=1\}\mbox{d}s\\
&=\lim_{t\to\infty}\frac{1}{t}\mathbb{E}\int_0^t\mathbf{I}\{Q_1(s-)={n}-1,
Q_2(s-)=1\}\mbox{d}A(s)\\
&\ \ \ \ +
\lim_{t\to\infty}\frac{1}{t}\mathbb{E}\int_0^t\mathbf{I}\{Q_1(s-)={n},
Q_2(s-)=0\}\mbox{d}A(s).
\end{aligned}
\end{equation}
\end{thm}

\begin{proof}
The proof of this theorem easily follows from representations
\eqref{LDF13}, \eqref{LDF14} and \eqref{LDF15}. We divide the
left- and right-hand sides by $t$, and pass to the limits in
probability as $t$ increases to infinity.

Specifically, \eqref{LDF16} and \eqref{LDF17} follow from
\eqref{LDF13} and from the fact that the martingales $m_{i,j}(t)$
vanish with probability 1 as $t\to\infty$. The last is the
consequence of vanishing $m_{\Pi_i}(t)$ and $m_{\Pi_{i+1}}(t)$
with probability 1 as $t\to\infty$. (The right-hand side of
\eqref{LDF13A} is divided by $t$, and then we pass to the limit
with probability 1 as $t\to\infty$.) Relations \eqref{LDF18} and
\eqref{LDF19} follow similarly from the corresponding relations
\eqref{LDF14} and \eqref{LDF15} and vanishing the martingales
$m_{{n},0}(t)$ and $m_{{n},1}(t)$ given by \eqref{LDF14A} and
\eqref{LDF15A}. We also use the Lebesgue dominated convergence
theorem to replace $\mathbb{P}^{\_}\lim_{t\to\infty}$ by
$\lim_{t\to\infty}\mathbb{E}$ in the places where it is required.
\end{proof}

The standard particular cases from Theorem \ref{thm1} are the
cases where $A(t)$ is Poisson, $D(t)$ is Poisson, and both $A(t)$
and $D(t)$ are Poisson all together. We do not consider the first
two cases, and the only last of these three cases is formulated
and proved below. In this case the representation is deduced
explicitly and has the form of limiting distributions rather than
frequencies.

\begin{cor}\label{cor1}
In the case where $A(t)$ is Poisson with rate $\lambda$, and
$D(t)$ is Poisson with rate $\Delta$ denote
$P_{i,j}=\lim_{t\to\infty}\mathbb{P}\{Q_1(t)=i, Q_2(t)=j\}$. We
have the following.

In the boundary cases $i=0$ for $j=0,1$ we have:
\begin{equation}
\label{LDF20} \mu P_{1,j}=\lambda P_{0,j}+j\Delta P_{0,j}.
\end{equation}

In the cases $i=1,2,\ldots,{n}-1$, $j=0,1$ we have:
\begin{equation}
\label{LDF21} i\mu P_{i,j}-(i+1)\mu P_{i+1,j}=\lambda
P_{i-1,j}-\lambda P_{i,j}+(1-j)\Delta P_{i-1,1}-jP_{i,1}.
\end{equation}

In the case $i={n}$ and $j=0$ we have:
\begin{equation}
\label{LDF22} {n}\mu P_{n,0} = \lambda P_{{n}-1,0} - \lambda
P_{{n},0} +\Delta P_{{n},1}.
\end{equation}

In the case $i={n}$ and $j=1$ we have:
\begin{equation}
\label{LDF23} {n}\mu P_{{n},1}=\lambda P_{{n}-1,1}+\lambda
P_{{n},0}.
\end{equation}
\end{cor}

\begin{proof}
In the case where $A(t)$ and $D(t)$ are Poisson, the proof of this
corollary is based on the following two relations:
\begin{equation}
\label{LDF24}
\begin{aligned}
&\lim_{t\to\infty}\frac{1}{t}\mathbb{E}\int_0^t\mathbf{I}\{Q_1(s-)=i,
Q_2(s-)=j\}\mbox{d}A(s)=\lambda P_{i,j},\\
&\ \ \ \ i=0,1,\ldots,{n}; \ j=0,1,
\end{aligned}
\end{equation}
and
\begin{equation}
\label{LDF25}
\begin{aligned}
&\lim_{t\to\infty}\frac{1}{t}\mathbb{E}\int_0^t\mathbf{I}\{Q_1(s-)=i,
Q_2(s-)=j\}\mbox{d}D(s)=\Delta P_{i,j},\\
&\ \ \ \ i=0,1,\ldots,{n}; \ j=0,1.
\end{aligned}
\end{equation}
Then inserting \eqref{LDF24} and \eqref{LDF25} into
\eqref{LDF16}-\eqref{LDF19} leads to the statement of the
corollary. Thus there is only required to prove \eqref{LDF24} and
\eqref{LDF25}.

Prove \eqref{LDF24}. By semimartingale decomposition of Poisson
process we have $A(t)=\lambda t+M_A(t)$. Substituting this
representation for \eqref{LDF24} we obtain:
\begin{equation}
\label{LDF26}
\begin{aligned}
&\lim_{t\to\infty}\frac{1}{t}\mathbb{E}\int_0^t\mathbf{I}\{Q_1(s-)=i,
Q_2(s-)=j\}\mbox{d}A(s)\\
&=\lim_{t\to\infty}\frac{\lambda}{t}\int_0^t\mathbb{P}\{Q_1(s-)=i,
Q_2(s-)=j\}\mbox{d}s\\
&\ \ \ \
+\lim_{t\to\infty}\frac{1}{t}\mathbb{E}\int_0^t\mathbf{I}\{Q_1(s-)=i,
Q_2(s-)=j\}\mbox{d}M_A(s)\\
&\ \ \ \ i=0,1,\ldots,{n}; \ j=0,1.
\end{aligned}
\end{equation}
The first term of the right-hand side is equal to $\lambda
P_{i,j}$, and the second term is equal to 0. Therefore
\eqref{LDF24} is proved.

The proof of \eqref{LDF25} is analogous. The corollary is proved.
\end{proof}

\section{Performance analysis and algorithms}\label{Performance}

One of the most important performance characteristics is the loss
proportion per arriving customer. Specifically, denoting the loss
proportion by $f$, we write
\begin{equation}
\label{PAA1}f=
\frac{\mathbb{P}^{\_}\lim_{t\to\infty}q_3(t)}{\mathbb{P}^{\_}\lim_{t\to\infty}a(t)}
=\frac{1}{\lambda}\ \mathbb{P}^{\_}\lim_{t\to\infty}q_3(t).
\end{equation}
The parameters $\lambda$, $\mu$ and $\Delta$ are assumed to be
given. So, the problem is to find the appropriate number of
servers $n$ such that $f\leq \alpha$. More accurately, the problem
is to minimize ${n}$ subject to $f\leq\alpha$.

In general case, we have no explicit representation for the
processes, and then simulation techniques are used. There are
three equivalent relations obtained for
$\mathbb{P}^{\_}\lim_{t\to\infty}q_3(t)$ in Section \ref{Basic
equations}: \eqref{SDN8}, \eqref{SDN9} and \eqref{SDN10}.

The simplest of them looks \eqref{SDN10}. According to
\eqref{SDN10} the appropriate limit
$\mathbb{P}^{\_}\lim_{t\to\infty}q_3(t)$ can be obtained
straightforwardly. Immediately before each arrival (one occurring
say at instant $s$) we only check the event $\{Q_1(s-)={n}\ \cap \
Q_2(s-)=1\}$. Such type of simulation need not require any theory.
However, in the case of small $\alpha$ this simulation is not
effective. As the large number of servers is chosen, and the loss
proportion should be small, resulting in erroneous conclusion.

The other relation is given by \eqref{SDN9}. For a small $\alpha$
it is also based on the rare events, and its application is
impossible.

In this general case only \eqref{SDN8} is available. It is based
on the average number of busy servers $\lim_{t\to\infty}
{t}^{-1}\mathbb{E}\int_0^tQ_1(s)\mbox{d}s$ which can be estimated
by simulation. More specifically, the simulation is based on
application of Theorem \ref{thm1}, helping us to obtain the
frequencies
\begin{equation}\label{PAA2}
\begin{aligned}
\lim_{t\to\infty}\frac{1}{t}\int_0^t\mathbb{P}\{Q_1(s)=i\}\mbox{d}s=&
\lim_{t\to\infty}\frac{1}{t}\int_0^t\mathbb{P}\{Q_1(s)=i,
Q_2(s)=0\}\mbox{d}s\\
&+\lim_{t\to\infty}\frac{1}{t}\int_0^t\mathbb{P}\{Q_1(s)=i,
Q_2(s)=1\}\mbox{d}s,\\
&\ \ \ \ i=0,1,\ldots,{n}.
\end{aligned}
\end{equation}
The first and second terms of the right-hand size are obtained
from Theorem \ref{thm1}. They are expressed via the frequencies
\begin{equation*}
\begin{aligned}
\label{PAA3}&\lim_{t\to\infty}\frac{1}{t}\mathbb{E}\int_0^t\mathbf{I}
\{Q_1(s-)=i, Q_2(s-)=j\}\mbox{d}A(s)\\
&\lim_{t\to\infty}\frac{1}{t}\mathbb{E}\int_0^t\mathbf{I}
\{Q_1(s-)=i, Q_2(s-)=j\}\mbox{d}D(s)\\
&i=0,1,\ldots,{n}; \ j=0,1,
\end{aligned}
\end{equation*}
which in turn are easily evaluated by simulation.

It is worth noting that direct simulating of
$\lim_{t\to\infty}t^{-1}\int_0^t\mathbb{P}\{Q_1(s)=i,
Q_2(s)=j\}\mbox{d}s$ ($i=0,1,\ldots,{n}$, $j=0,1$) and
consequently $\lim_{t\to\infty}
{t}^{-1}\mathbb{E}\int_0^tQ_1(s)\mbox{d}s$ is not available. (See
Section 9 of \cite{Abramov 2006} for the detailed discussion of
this question.)

Thus, for the simulation purposes we use \eqref{SDN8} and the
relations of Theorem \ref{thm1}.

Recall that the problem is to minimize $m$ subject to $f\leq
\alpha$. The upper bound for ${n}$ is unknown, and the special
search procedure is necessary. The relevant search procedure is
known due to Rubalskii \cite{Rubalskii 1982}, who proposed the
search algorithm for minimization of a unimodal function on an
unbounded set. The optimal algorithm is an extension of the
standard Fibonacci procedure.

Thus, solution of the problem is based on two main steps:

\smallskip $\bullet$ Simulation

\smallskip $\bullet$ Search step

\smallskip

These procedures are repeated until the optimal solution is not
found.

\smallskip
In the particular case where both $A(t)$ and $D(t)$ are Poisson,
for stationary probabilities $P_{i,j}$, $i=0,1,\ldots,{n}$;
$j=0,1$, we have the system of algebraic equations. Furthermore,
the upper bound for the value ${n}$ can be evaluated as well. In
order to evaluate the upper bound for ${n}$, one can imagine that
arrival rate is $\lambda+\Delta$, that is retrials occur all the
time continuously and the retrial space is always occupied. In
this case the loss probability is greater then that original, and
consequently the corresponding value ${n}$ is greater than that
required. This upper bound ${n}$  can be easily calculated by the
Erlang loss formula. Then, having the upper and lower bounds we
obtain the optimal number of servers easier.

The similar approach can be used and for non-Markovian system. In
the most practical cases where input stream is recurrent the upper
bound for $n$ can be easily evaluated by known formulae for the
loss systems.

However in general we have to use one or other special algorithm,
say the method of \cite{Rubalskii 1982}, in order to find optimal
value of $n$.

\section{Numerical examples}\label{Numerical examples}
In this section we provide numerical solution for two systems. The
first system is the Markovian system with $\lambda=10$,
$\Delta=2$, $\mu=1$ and $\alpha=0.0001$ In the second example we
assume that the increments of the processes $A(t)$ and $D(t)$ are
deterministic with the same parameters.

\subsection{Example 1} Assuming that the lower bound for $m$ is
equal to 1, let us find the upper bound. For this purpose we
consider the M/M/${n}$/0 queueing system with input rate
12=$\lambda$+$\Delta$ and $\mu=1$. According to Erlang loss
formula, for the loss probability we have
\begin{equation*}
\label{NE1}p_{n}=\dfrac{\dfrac{\rho^{n}}{{n}!}}{\sum_{i=0}^{n}\dfrac{\rho^i}{i!}},
\end{equation*}
where $\rho=\lambda/\mu$=12. Assuming that the loss probability is
equal to $\alpha=0.0001$ let us find the value ${n}$.

For ${n}$=26 the loss probability is 0.000174, and for ${n}$=27
the loss probability is 0.000078 The first of these probabilities
is greater than 0.0001, and the second one is smaller than 0.0001
Therefore, the upper bound for ${n}$ is 27.

Let us now find ${n}$ optimal. The solution of this problem can be
achieved as follows. We find a new value of ${n}$=1+(27-1)/2=14.
Now we solve the system of equations \eqref{LDF20}-\eqref{LDF23}
given by Corollary \ref{cor1}. We also use \eqref{PAA1} and
\eqref{SDN13}. For loss probability we obtain 0.064841 Therefore
in the next step we find a new value ${n}$ by 14+(27-14)/2=20.5
taking then the integer part. With new value ${n}=20$ we solve the
system again, and for the loss probability we obtain 0.001226 This
value is greater than $\alpha$, and therefore the new value of
${n}$ is the integer part of 20+(27-20)/2 equal to 23. With this
value of ${n}$ for the loss probability we have 0.000110 Thus the
last step is anticipated with ${n}$=24. (Formally we must analyze
first the value ${n}$=25 and only then accept ${n}$=24.) For this
value ${n}$ the loss probability is 0.000045 Thus the problem is
solved, and the value ${n}$=24 is the optimal number of servers
decreasing the loss probability to the required value.

\subsection{Example 2} Assume that all parameters are the same as
in Example 1, but increments of the processes $A(t)$ and $D(t)$
are deterministic. In this case the scheme of calculation includes
simulation as it is explained in Section \ref{Performance}. For
this concrete system the algorithm can be simplified.
Specifically, as in Example 1 we first evaluate the upper bound
for ${n}$. For our approximation of the upper bound of ${n}$ we
consider the D/M/${n}$/0 queueing system, the interarrival times
of which all are equal to 1/($\lambda+\Delta$), i.e. in our case
1/12. Notice, that under the assumption that customers constantly
arrive to the main queue by $A(t)$ and $D(t)$ we do not have the
model of D/M/${n}$/0 queue exactly. Although input stream is based
on two sources with deterministic interarrival times, the
structure of the overall input stream is complicated. The
arguments for such approximating by the D/M/${n}$/0 queueing
system are heuristic. Nevertheless, such approximation is
available. This is supported by computations as well.

The loss probability for GI/M/${n}$/0 loss systems is well-known
(e.g. \cite{Cohen 1957}, \cite{Palm 1943}, \cite{Pollaczek 1953},
\cite{Takacs 1957} as well as \cite{Bharucha-Reid 1960}):
\begin{equation}\label{NE2}
p_{{n}}= \left[\sum_{i=0}^{n}\binom{{n}}{i} \prod_{j=1}^i
\frac{1-r_{j}}{r_{j}}\right]^{-1},
\end{equation}
where
\begin{equation*}
r_{j}=\int_0^\infty\mbox{e}^{-j\mu x}\mbox{d}R(x),
\end{equation*}
and $R(x)$ is the probability distribution function of an
interarrival time. In our case interarrival times are
deterministic of length 1/12 and $\mu=1$. Therefore, in
\eqref{NE2}
\begin{equation*}
r_j=\mbox{e}^{-j/12}.
\end{equation*}

According to our calculations the upper bound for ${n}$ is 22. The
value of the loss probability is 0.000036 Let us note that for
${n}$=21 the value of the loss probability is 0.000137, i.e. it is
slightly greater than $\alpha$=0.0001

In the next step we choose then the integer part of 1+(22-1)/2,
and the value ${n}$ is 11. Simulating with this value of ${n}$
yields the value of the loss probability 0.125444 Therefore in the
next step we choose the integer part of 11+(22-11)/2, i.e.
${n}$=16. With this value of ${n}$ by simulating we have the value
of loss probability 0.002784 In the next step
${n}$=16+(22-16)/2=19. With this value ${n}$ we correspondingly
obtain the loss probability 0.000060 In the next step we choose
the integer part of 19-(19-16)/2 for ${n}$, i.e. ${n}$=17. The
loss probability in this case is 0.000786 There is one value
${n}$=18 which must be checked. For this value the loss
probability is 0.000214 Thus our conclusion is ${n}$=17. Recall
that the loss probability in this case is 0.000060

\section{Concluding remarks}\label{Concluding remarks}
In this paper we studied multiserver retrial queueing system
having only one space in orbit. For this queueing system we
derived analytic results permitting us to provide performance
analysis. Specifically we solved the problem to find the possibly
minimal number of servers decreasing the loss probability to given
small value. In general case the algorithm of solution requires
effective simulation, meeting the challenge of rare events
problem. Numerical examples were provided for two cases. In the
first case all processes were assumed to be Poisson, and the
system Markovian. For Markovian system the calculations were based
on analytic results. In the second case the processes $A(t)$ and
$D(t)$ had deterministic increments. The analysis in this case was
based on effective simulation. Comparison of these results showed
that in the case of Markovian system the required number of
servers, making the loss proportions smaller than the given value,
is relatively greater than that in the system having the processes
$A(t)$ and $D(t)$ with deterministic increments.

% ----------------------------------------------------------------
%\bibliographystyle{amsplain}
%\bibliography{}

\end{document}